\documentclass[12pt]{amsart}
\usepackage{amsmath,amssymb,amsthm,amscd}
\numberwithin{equation}{section}

\def\SS{\Bbb S}

\newtheorem{prop}{Proposition}[section]
\newtheorem{theo}[prop]{Theorem}
\newtheorem{lemm}[prop]{Lemma}
\newtheorem{coro}[prop]{Corollary}
\newtheorem{rema}[prop]{Remark}
\newtheorem{exam}[prop]{Example}

\def\begeq{\begin{equation}}
\def\endeq{\end{equation}}
\def\p{\partial}

\def\lf{\left}
\def\ri{\right}

\def\ol{\overline}

\def\hbk{\mathbb{H}^{ 3}_{-\kappa^2 }}
\def\hb{\mathbb{H}^3_{-1}}
\def\ii{\sqrt{-1}}
\def\la{\langle}
\def\ra{\rangle}
\def\Minkowski{\mathbb{R}^{3,1}}
\begin{document}

\title{Boundary behaviors  and scalar curvature of compact manifolds}
\author{Yuguang Shi$^1$}
\thanks{$^1$Research partially supported by NSF grant of China and and
Fok YingTong Education Foundation.}

\address{Key Laboratory of Pure and Applied mathematics, School of Mathematics Science, Peking University,
Beijing, 100871, P.R. China.} \email{ygshi@math.pku.edu.cn}

\author{Luen-Fai Tam$^2$}
\thanks{$^2$Research partially supported by Earmarked Grant of Hong
Kong \#CUHK403005.}

\address{Department of Mathematics, The Chinese University of Hong Kong,
Shatin, Hong Kong, China.} \email{lftam@math.cuhk.edu.hk}
\renewcommand{\subjclassname}{%
  \textup{2000} Mathematics Subject Classification}
\date{October 2006}

\begin{abstract} In this paper, by modifying the arguments in
\cite{WY}, we get some  rigidity theorems on compact manifolds with
nonempty boundary. The results in this paper  are  similar  with
those in \cite{ST} and \cite{WY}. Like \cite{ST} and \cite{WY}, we
still use quasi-spherical metrics introduced by \cite{Ba} to get
monotonicity of some quantities.
\end{abstract}

\maketitle \markboth{Yuguang Shi and Luen-Fai Tam} {Boundary
behaviors  and scalar curvature of compact manifolds}

\section{Introduction}
In  \cite{ST}, the authors proved the following: Let $(\Omega, g)$
be a compact manifold of dimension three with smooth boundary
$\Sigma$ which has positive Gaussian curvature and has positive
mean curvature. Suppose $\Omega$ has nonnegative scalar curvature,
then for each boundary component $\Sigma_i$ of $\Sigma$ satisfies,
\begin{equation}\label{Brown-Yorkmass}
\int_{\Sigma_i}(H^i_0 -H)\,d\Sigma_i \ge 0
\end{equation}
where $H^i_0$ is the mean curvature of $\Sigma_i$ with respect to
the outward normal when it is isometrically embedded in $\mathbb{
R}^3$, $d\Sigma_i$ is the volume form on $\Sigma_i$ induced from
$g$. Moreover, if equality holds for some $\Sigma_i$ then $\Sigma$
has only one component and $\Omega$ is a domain in $\mathbb{R}^3$.

The result gives  restriction on a convex surface $\Sigma$ in
$\mathbb{R}^3$ which can bound a compact manifold with nonnegative
scalar curvature such that the mean curvature of $\Sigma$ is
positive. It is interesting to see what one can say for convex
surface $\Sigma$ in $\hbk$, the hyperbolic space with constant
curvature $-\kappa^2$.

The  result   mentioned above   has  other interpretation. It
implies    the quasi-local mass introduced by Brown-York
\cite{BY1,BY2} is positive under the condition that the boundary
has positive Gaussian curvature. In \cite{LY1}, \cite{LY2}, Liu
and Yau introduced  a quasi-local mass in spacetime.
 This quasi-local
mass were also introduced by Epp \cite{E} and Kijoswki \cite{K}.
More importantly, Liu and Yau  proved   its positivity, using
\cite{ST}. A recent definition of quasi-local mass that relates with
these works please see \cite{Z}.

Motivated by   \cite{ST,LY1,LY2}, in a recent work \cite{WY}  Wang
and Yau proves the following: Suppose $(\Omega,g)$ is a three
dimensional manifold with smooth boundary $\Sigma$ with positive
mean curvature $H$, which is a topological sphere. Suppose the
scalar curvature $\mathcal{R}$ of $\Omega$ satisfies
$\mathcal{R}\ge -6\kappa^2$ and the Gaussian curvature of $\Sigma$
is larger than $-\kappa^2$, then there is a future directed
time-like vector value function $\mathbf{W^0}$  on $\Sigma$ such
that
$$
\int_\Sigma (H_0-H)\mathbf{W}^0 d\Sigma
$$
is time-like. Here $H_0$ is the mean curvature of $\Sigma$ when
isometrically embedded in $\hbk$, which is in turns isometrically
embedded in $\Minkowski$, the Minkowski space. In this result, the
vector $\mathbf{W}^0$   is not very  explicit  because it is
obtained by solving a backward parabolic equation by prescribing
data at infinity.

 In this work, by modifying the argument in \cite{WY}, we   get
 similar result
by replacing $\mathbf{W}^0$ by   $\mathbf{W}_{\Sigma_0}= (x_1,
x_2, x_3, \alpha t)$ for some $\alpha>1$ depending only on the
intrinsic geometry of $\Sigma$. Here $(x_1, x_2, x_3, t)$ is the
future directed unit normal vector of $\hbk$ in $\mathbb{
R}^{3,1}$. See Theorem \ref{positivity2} for a more precise
statement.  We believe that the same result should be true with
$\mathbf{W}_{\Sigma_0} =(x_1, x_2, x_3, t)$, but we cannot prove
it for the time being.

As a consequence, if $o$ is a point inside  of $\Sigma$ in $\hbk$,
then
$$
\int_\Sigma (H_0-H)\cosh \kappa r \,d\Sigma\ge0.
$$
where $r$ is the distance function from $o$ in $\hbk$. Moreover,
equality holds if and only if $(\Omega,g)$ is a domain in $\hbk$.
The results can be considered as   generalization of the results
in \cite{ST}. In fact, if we let $\kappa\to0$, we may obtain the
inequality (\ref{Brown-Yorkmass}).

The paper is organized as follows. In \S2, we list some facts that
we need, most of them are from \cite{WY}. In \S3, we prove our main
results. We will also give some examples that $\alpha$ in Theorem
\ref{positivity2} can be taken to be $1$ and also study some
properties of $\int_\Sigma (H_0-H)\cosh \kappa r \,d\Sigma$.

\section{Preliminary}

Most materials in this section are from Wang and Yau \cite{WY}.
Let $(\Omega,g)$ be a compact manifold with smooth boundary so
that $\Sigma=\p\Omega$ is a topologically sphere. Let $H$ be the
mean curvature with respect to the outward normal and $K$ be the
Gaussian curvature of $\Sigma$ and let $\mathcal{R}$ be the scalar
curvature of $\Omega$. In our convention, the mean curvature of
the unit sphere in $\mathbb{R}^3$ with respect to the outward
normal is 2. By \cite{P1,CW}, we have the following:

\begin{lemm}\label{embedding1}  Suppose the Gaussian curvature
$K$ of $\Sigma$ satisfies  $K> -\kappa^2$. Then $\Sigma$
can be isometrically embedded into the hyperbolic space
$\hbk$ with constant curvature $-\kappa^2$ as a convex
surface which bounds a convex domain $D$ in $\hbk$.
Moreover, the embedding is unique up to an isometry of
$\hbk$.
\end{lemm}
Since $\Sigma$ is a topological sphere, its image $\Sigma_0$ under
the embedding divides $\hbk$ into two components, the exterior and
the interior of $\Sigma_0$. Let $\mathbf{N}$ be the unit outward
normal of $\Sigma_0$. It is said to be convex if the second
fundamental form $h(X,Y)=-\la \nabla^\kappa_XY,\mathbf{N}\ra$ is
positive definite for $X$ and $Y$ tangent to $\Sigma_0$, where
$\nabla^\kappa$ is the covariant derivative of $\hbk$. The
interior $D$ of $\Sigma_0$ being convex means that $D$ is
geodesically convex.

The existence and uniqueness of the embedding were proved by
Pogorelov \cite{P1}. The convexity of $\Sigma_0$ and $D$ were
proved by do Carmo and Warner \cite{CW}.

Further identify $\hbk$  with
$$
\{(x_1,x_2,x_3,t)\in \mathbb{R}^{3,1}|\
x_1^2+x_2^2+x_3^2-t^2=-\frac1{\kappa^2},\ t>0\}
$$
where $\mathbb{R}^{3,1}$ is the Minkowski space consisting
of $\mathbf{X}=(x_1,x_2,x_3,t)$ with the Lorentz metric
$dx_1^2+dx_2^2+dx_3^2-dt^2$. Position vectors in
$\mathbb{R}^{3,1}$ can be parametrized as
\begin{equation}\label{positionvector}
\begin{split}\mathbf{X} & =(x_1,x_2,x_3,t)\\
&=\frac1\kappa(\sinh\kappa r\cos\theta,\sinh\kappa
r\sin\theta\cos\psi,\sinh\kappa r\sin\theta\sin\psi,\cosh
\kappa r).\end{split}
\end{equation}
The  metric of $\hbk$ is then $$dr^2+ \kappa^{-2} \sinh^2
\kappa r d\sigma^2=dr^2+ \kappa^{-2} \sinh^2 \kappa r
(d\theta^2+\sin^2\theta d\psi^2).$$ Note that $r$ is the
geodesic distance of a point from $(0,0,0,1/\kappa)\in
\hbk$.

Let $\Sigma_\rho$ be the level surface outside $\Sigma_0$ in
$\hbk$ with distance $\rho$ from $\Sigma_0$. Foliate
$\hbk\setminus D$ by $\Sigma_\rho$, $\rho\ge0$. The hyperbolic
metric can be written as $ds^2_{\hbk}=d\rho^2+g_\rho$, where
$g_\rho$ is the induced metric on $\Sigma_\rho$.

Suppose f $F:\Sigma\to \hbk$  is the embedding with unit outward
normal $\mathbf{N}$. Then $\Sigma_\rho$ as a subset of
$\mathbb{R}^{3,1}$ is given by
\begin{equation}\label{embedding2}
  \mathbf{X}(p,\rho)=\cosh\kappa\rho\,
  \mathbf{X}(p,0)+\kappa^{-1}\sinh\kappa\rho\,
  \mathbf{N}(p,0)
\end{equation}
Here for simplicity, $(p,\rho)$ denotes a point
$\Sigma_\rho$ which lies on the geodesic from the point
$p\in \Sigma_0$ and $\mathbf{X}(p,0)=\mathbf{X}(F(p))$.

Suppose in addition that the mean curvature of $H$ of $\Sigma$ with
respect to $(\Omega,g)$ is positive and the scalar curvature
$\mathcal{R}$ of $\Omega$ is greater than or equal to $-6\kappa^2$.
Wang and Yau \cite{WY} are able to solve the following parabolic
equation
\begin{equation}\label{parabolicequation}
\left\{
  \begin{array}{ll}
    2H_0 \frac{\p u}{\p \rho}=2u^2 \Delta_\rho u + (u-u^3)(R^\rho+6\kappa^2), \\
    u(p, 0)=\frac{H_0 (p, 0)}{H(p)},
  \end{array}
\right.
\end{equation}
for all $\rho\ge0$ with positive and bounded solution $u$. Here
$\Delta_\rho$ is the Laplacian operator of $\Sigma_\rho$, $R^\rho$
is scalar curvature of $\Sigma_\rho$, and $H_0$ is the mean
curvature of $\Sigma_\rho$ which is positive. We need the
following result which is proved   by Wang and Yau, see Theorem
6.1 and Corollary 6.1 in \cite{WY}.
\begin{theo}\label{positivemass}{\rm [Wang-Yau]} Let $(\Omega, g)$ be a $3$-dimensional  compact
 Riemannian manifold
with nonempty smooth boundary which is a topological sphere.
Suppose the scalar curvature of $\mathcal{R}$ of $\Omega$
satisfies $\mathcal{R}\ge -6\kappa^2$, the Gaussian curvature of
its boundary $\Sigma$ satisfies $K> -\kappa^2$, and the mean
curvature $H$ of the boundary with respect to outward unit norm is
positive. Let $\mathbf{X}$ be the position vector of
$\mathbf{H}^3_{-\kappa^2}$ in $\mathbb{ R}^{3,1}$ then
$$\lim_{\rho\rightarrow\infty}\int_{\Sigma_\rho}(H_0
-H)\mathbf{X}\cdot\zeta\le0
$$
 for any  future directed null vector $\zeta$ in $\mathbb{
 R}^{3,1}$.
\end{theo}

\begin{coro}\label{positivemass1} With the same assumptions and
notations as in Theorem \ref{positivemass},
$$\lim_{\rho\rightarrow\infty}\int_{\Sigma_\rho}(H_0
-H)\cosh\kappa r d\Sigma_\rho\ge0.
$$
where $r$ is as in (\ref{positionvector}).
\end{coro}

We also have the following rigidity result.

\begin{prop}\label{rigidity}  With the same assumptions and notations as in
Theorem \ref{positivemass}. Suppose
\begin{equation}\label{rigidity1}
\lim_{\rho\to\infty}\int_{\Sigma_\rho}(H_0
-H)\mathbf{X}\cdot\zeta\,d\Sigma_\rho=0
\end{equation}
for some  future directed null  vector $ \zeta$ in
$\mathbb{R}^{3,1}$, where the inner product is given by the
Lorentz metric. Then $\Omega$ is a domain in $\hbk$.
 \end{prop}
 \begin{proof} For simplicity, let us assume that $\kappa=1$.
As in the proof of Theorem \ref{positivemass} in \cite{WY}, let
$ds^2=u^2d\rho^2+g_\rho$ be the quasi-spherical metric where $u$
is the solution of (\ref{parabolicequation}). Let $(M,\tilde g)$
be the manifold by gluing $(\Omega,g)$ with $\hb\setminus D$ with
metric $ds^2$.
 By
 \cite{WY}, if (\ref{rigidity1}) is true, then the manifold $(M,\tilde g)$
 has a Killing spinor $\phi$ which is nontrivial, smooth away from
 $\Sigma$ and is continuous. More precisely, $\phi$ satisfies:

 \begin{equation}\label{spinoreq}
 \nabla_V\phi +\frac{\ii}2c(V)\cdot\phi=0
\end{equation}
 where $c(V)\cdot$ is the Clifford multiplication.
 Hence $M\setminus\Omega$ is Einstein.
 Since $M$ has dimension three,  the sectional curvature is $-1$ in
 $(M\setminus\Omega,\tilde g)$, see \cite{B}, for example.  Let $h^0_{ij}$ and $h_{ij}$ be the
 second fundamental form of $\Sigma_0$ with respect to the metrics
 $ds_{\hb}^2$ and $ds^2$ respectively. Then
 $h_{ij}=u^{-1}h_{ij}^0$. By the Gauss equation and the fact that
 both $ds_{\hb}^2$ and $ds^2$ have constant curvature $-1$,
 $u\equiv1$.

 On the other hand, $\phi$ is not zero on $\Sigma_0$
 and so $\phi$ is a nontrivial Killing spinor in $(\Omega,g)$
 satisfying (\ref{spinoreq})
 and   $g$ has constant curvature $-1$ as before.

  We claim that the second fundamental
 forms of $\Sigma_0$ with respect to $g$ and $ds_{\hb}^2$ are
 equal. If this is true, then by the proof of \cite[Lemma 4.1]{ST}, we can conclude that $\tilde g$ is smooth. From this
 it is easy to see that $(M,\tilde g)=\hb$.

 Let us prove the claim. Since $\Sigma_0$ is a topological sphere,
  for some $a>0$
 a tubular neighborhood of $\Sigma_0\times(0,a)$ in $(\Omega,g)$ is simply
 connected. Since it has constant curvature $-1$,
 $(\Sigma_0\times(0,a),g)$ can be isometrically embedded in $\hb$,
 see \cite[p.43]{Sp}. Denote the embedding by
 $\iota=(u_1,u_2,u_3)$ where $(u_1,u_2,u_3)$ is are global
 coordinates in $\hb$. Since $\iota$ is an isometry, the normal
 curvatures of $\Sigma_0\times\{\tau\}$ for $0<\tau<a$ with respect
 to $g$ and the hyperbolic metric are equal. Hence they are
 uniformly bounded on $\Sigma_0\times(0,a)$. Note that $\Sigma_0$
 is convex in $\hb$, so $\Sigma_0\times\{\tau\}$ is also convex
 when embedded in $\hb$, for $0<\tau<a$ provided $a$ is small. By \cite[VI, \S3]{P2}, for any $k\ge0$,
 $|\nabla_{\tau}^ku_i|$ are uniformly bounded on
 $\Sigma_0\times(0,a)$, where $\nabla_\tau$ is the covariant
 derivatives of $\Sigma_0\times\{\tau\}$ with induced metric by
 $g$. Hence by taking a subsequence of $\tau_j\to0$, we obtain an
 isometric embedding of $(\Sigma,g)$. In this embedding the
 second fundamental form with respect to $g$ and $ds^2_{\hb}$ are
 equal. Since the embedding of $\Sigma$ is unique up to an isometry of $\hb$,
 the claim is true.
 \end{proof}

\section{Main results}
Let $(\Omega,g)$ be as in Theorem \ref{positivemass} and
let $\p\Omega=\Sigma$. With the same notations as in the
theorem, suppose $o=(0,0,0,1/\kappa)$ is in $D$ which is
the interior of $\Sigma$ in $\hbk\subset \Minkowski$.   Let
$\Sigma_0$ be the image of $\Sigma$   under the embedding
described in \S1.  $H_0$ is the mean curvature of
$\Sigma_0$  in $\hbk$. We also identify $\Sigma$ with
$\Sigma_0$. Let $B_o(R_1)$ and $B_o(R_2)$ be geodesic balls
in $\hbk$ such that $B_o(R_1)\subset D\subset B_o(R_2)$.

 We want to prove the following:
 \begin{theo}\label{positivity1}
Let $(\Omega, g)$ be a compact manifold with smooth
boundary $\Sigma$. Assume the following are true
\begin{itemize}
  \item[(i)] The scalar curvature $\mathcal{R}$ of $(\Omega,g)$ satisfies
  $\mathcal{R}\ge -6\kappa^2$ for some $\kappa>0$.
  \item[(ii)] $\Sigma$ is a topological sphere with
  Gaussian curvature $K>-\kappa^2$ and with positive mean
  curvature $H$.
\end{itemize}
 With the above notations, for any future directed null vector
 $\zeta$ in $\Minkowski$,
 \begin{equation}\label{positivity2}
m(\Omega;\zeta)=\int_{\Sigma}(H_0-H)\mathbf{W}_{ \Sigma_0}\cdot
\zeta d\Sigma\le0
 \end{equation}
where $\mathbf{W}_{\Sigma_0}=(x_1,x_2,x_3,\alpha t)$ with
$$
\alpha=\coth \kappa R_1+\frac1{\sinh \kappa R_1} \lf(\frac{\sinh^2
\kappa R_2}{\sinh^2 \kappa  R_1}-1\ri)^\frac12,
$$
$\mathbf{X}=(x_1,x_2,x_3, t)$ being the position vector in
$\Minkowski$ and the inner product is given by the Lorentz metric.
Moreover, if equality holds in (\ref{positivity2}) for some future
directed null vector $\zeta$, then $(\Omega,g)$ is a domain in
$\hbk$.
\end{theo}

\begin{rema}   If $R_1, R_2\to\infty$ in such a way that
$R_2- 2R_1\to -\infty$, then $\alpha\to1$.
\end{rema}

In (\ref{positionvector}), the position vector in $\Minkowski$ is
given by \begin{equation}\label{positionvector1}  \mathbf{X}
=(x_1,x_2,x_3,t)
 =\frac1\kappa(\phi_1\sinh\kappa r,\phi_2\sinh\kappa
r,\phi_3\sinh\kappa r,\cosh \kappa r)
\end{equation}where $(\phi_1,\phi_2,\phi_3)$
denote  position vectors of points of $\mathbb{S}^2$ in
$\mathbb{R}^3$. Let $\{\Sigma_\rho\}$ be the foliation of
$\hbk\setminus D$ described in \S1. We need the following:

\begin{lemm}\label{estimates} With the same assumptions and notations as in
Theorem \ref{positivity2}, let $y_1, y_2, y_3\in
\mathbb{R}$ with $\sum_{i=1}^3y_i^2=1$, the following are
true.
\begin{itemize}
  \item [(i)] For any $\rho>0$,
  \begin{equation}\label{monotone1}
   \frac{\p r}{\p\rho}\ge \frac{\sinh \kappa R_1}{\sinh \kappa R_2}
\end{equation}
  \item [(ii)] If $\phi=\sum_{i=1}y_i\phi_i$, then for
  $\rho>0$
   \begin{equation}\label{monotone2}
(\frac{\p\phi}{\p\rho})^2\le(1-\phi^2)\kappa^2\sinh^{-2}\kappa
r\lf(1-(\frac{\p r}{\p\rho})^2\ri)
\end{equation}
\end{itemize}
Hence
\begin{equation}\label{monotone2'}
\mu\cdot\kappa\frac{\p r}{\p\rho}\ge  \lf|\frac{\p\phi}{\p\rho} \ri|
\end{equation}
where
$$
\mu=\frac1{\sinh \kappa R_1} \lf(\frac{\sinh^2 \kappa R_2}{\sinh^2
\kappa  R_1}-1\ri)^\frac12.
$$
\end{lemm}
\begin{proof}
(i) Recall that   $o=(0,0,0,1/\kappa)\in D$ and $r$ is the
geodesic distance in $\hbk$ from $o$.   $D$ is geodesically
convex by Lemma \ref{embedding1}. For any $x, y\in \hbk$,
denote the geodesic from $x$ to $y$ parametrized by arc
length by $\ol{xy}$.

 Let $p\in
\Sigma$ and let $\gamma(\rho)$ be the geodesic through $p$ so that
$\gamma(\rho)$ is orthogonal to $\Sigma$ at $p$ with arc length
parametrization. Moreover, $\gamma(0)=p$ and $\gamma(\rho)$ is
outside $\Sigma$ for $\rho>0$. Let $q$ be the point on $\gamma$
such that $d(o,q)=d(o,\gamma)$. Then $q=\gamma(\rho_1)$ with
$\rho_1<0$ because $D$ is geodesically convex. Since $\frac{\p
r}{\p \rho} =\langle\nabla r, \nabla \rho\rangle$, $\frac{\p
r}{\p\rho}$ is nondecreasing on $\rho>0$ along $\gamma$ and
$\frac{\p r}{\p\rho}>0$ at $p$.   Let $\beta=\ol{op}$  and let $
\eta$ be the geodesic from $p$ on the totally geodesic
$\mathbb{H}^2_{ -\kappa^2}$ containing $\gamma$ and $\beta$ such
that $\eta$ is tangent to $\Sigma_0$. Let $x, y$ be the
intersection of $\eta$ with $\p B_o(R_2)$.  Then $\frac{\p
r}{\p\rho}\ge \sin \varphi$ where $ \varphi$ is the angle between
$\ol{xo}$ and $\ol{xp}$. Since $\eta$ is outside $B_o(R_1)$ so
$$
\sin\varphi\ge\frac{\sinh \kappa R_1}{\sinh \kappa R_2}.
$$
Hence we have
\begin{equation}\label{monotone3}
   \frac{\p r}{\p\rho}\ge \frac{\sinh \kappa R_1}{\sinh \kappa R_2}
\end{equation}
on $\hbk\setminus D$.

(ii) Since the inner product in $ \mathbb{R}^3$ of $(y_1,y_2,y_3)$
and $(\phi_1,\phi_2,\phi_3)$ is $\phi$, we may assume that
$\phi=\cos\theta$ in (\ref{positionvector}). The hyperbolic metric
outside $D$ is given by
$$
ds^2_{\hbk}=d\rho^2+g_\rho=dr^2+\frac1{\kappa^2}\sinh^2\kappa
r(d\theta^2+\sin^2\theta d\psi^2).
$$
Compute $ds^2_{\hbk}(\frac{\p}{\p\rho},\frac{\p}{\p\rho})$ in the
above two forms of $ds^2_{\hbk}$, we have
$$
 1=(\frac{\p r}{\p\rho})^2+\frac1{\kappa^2}\sinh^2 \kappa
r\lf[(\frac{\p\theta}{\p\rho})^2+\sin^2\theta(\frac{\p\psi}{\p\rho})^2\ri]
\ge
(\frac{\p r}{\p\rho})^2+\frac1{\kappa^2}\sinh^2 \kappa r
(\frac{\p\theta}{\p\rho})^2.
 $$
Since $\phi=\cos\theta$, (ii) follows.

The last assertion follows from (i), (ii), the fact that
$|\phi|\le1$ and the fact that $r\ge R_1$ for $\rho\ge 0$.

\end{proof}
\begin{lemm}\label{coordinatefunction1} With the same assumptions
and notations as in Theorem \ref{positivity2},
$$
H_0\frac{\p}{\p \rho}\mathbf{X}+
\Delta_{\rho}\mathbf{X}-2\kappa^2\mathbf{X}=\mathbf{0}
$$
in $\hbk\setminus D$.
\end{lemm}
\begin{proof} In the representation of $\hbk$ in
(\ref{positionvector1}), the Laplacian of $\hbk$ is given by
$$\Delta=\frac{\p^2}{\p r^2}+2\kappa\coth\kappa r\frac{\p}{\p r}+
\kappa^2\sinh^{-2}\kappa r \Delta_{\mathbb{S}^2}.
$$
So $\Delta \mathbf{X}=3\kappa^2\mathbf{X}$. In the  foliation
(\ref{embedding2}),  the metric of $\hbk$ is given by
$d\rho^2+g_\rho$ where $g_\rho$ is the induced metric on level
surface $ \Sigma_\rho$. The Laplacian on $\hbk$ is given by
$$
\Delta=\frac{\p^2}{\p \rho^2}+H_0\frac{\p}{\p \rho}+
\Delta_{\rho}.
$$
Using (\ref{embedding2}), we have
\begin{equation}\label{coordinatefunction}
\begin{split}
3\kappa^2\mathbf{X}&=\frac{\p^2}{\p
\rho^2}\mathbf{X}+H_0\frac{\p}{\p
\rho}\mathbf{X}+ \Delta_{\rho}\mathbf{X}\\
&=\kappa^2\mathbf{X}+H_0\frac{\p}{\p \rho}\mathbf{X}+
\Delta_{\rho}\mathbf{X}.
\end{split}
\end{equation}
From this the result follows.
\end{proof}
Now we are ready to prove Theorem \ref{positivity2}.

\begin{proof}[Proof of Theorem \ref{positivity2}]
 By \cite{WY} and Lemma \ref{coordinatefunction1}
\begin{equation}\label{monotone4}
    \begin{split}\frac{d}{d\rho}&\int_{\Sigma_\rho}(H_0-H)\mathbf{X}d
\Sigma_\rho\\&=\int_{\Sigma_\rho}-\frac12u^{-1}(u-1)^2(R^\rho+2\kappa^2)\mathbf{X}+
(u-1)\lf(\frac{H_0}u \frac{\p}{\p \rho}\mathbf{X}+\Delta_\rho
\mathbf{X}-2\kappa^2\mathbf{X}\ri)d
\Sigma_\rho\\
&=-\int_{\Sigma_\rho}u^{-1}(u-1)^2\lf[\frac12(R^\rho+2\kappa^2)\mathbf{X}
+H_0\frac{\p}{\p \rho}\mathbf{X}\ri]d \Sigma_\rho.
\end{split}
\end{equation}
Let $\lambda_a(p,\rho)$ be the principal curvature of the level
surface. Then $\lambda_a=\kappa\tanh \kappa(\mu_a+\rho)$, $\kappa$
or $\kappa\coth \kappa(\mu_a+\rho)$ with $\mu_a>0$, see \cite{WY}.
Hence $H_0=\lambda_1+\lambda_2$ and $R^\rho+2\kappa^2=2 \lambda_1
\lambda_2$. Let $W_0=\alpha\cosh \kappa r$ and $W=\phi\sinh\kappa
r$ where $\phi=\sum_{i=1}^3y_i\phi_i$, $\sum_{i=1}^3y_i^2=1$ and
$\alpha$ is the constant in the theorem. Then
\begin{equation}\label{monotone5}
\frac12(R^\rho+2\kappa^2)W_0 +H_0(W_0)_\rho =\alpha\lf[
 \lambda_1 \lambda_2\cosh\kappa r+\kappa(\lambda_1
+\lambda_2)\sinh\kappa r\cdot\frac{\p r}{\p \rho}\ri],
\end{equation}
\begin{equation}\label{monotone6}
\begin{split}
\frac12(R^\rho+2\kappa^2)W +H_0W_\rho &= \lambda_1
\lambda_2\phi\sinh\kappa r+\\
&\quad  \kappa(\lambda_1 +\lambda_2)\lf(\phi\cosh\kappa r
\cdot\frac{\p r}{\p \rho}+\frac1\kappa\sinh\kappa
r\cdot\frac{\p\phi}{\p \rho}\ri).
\end{split}
\end{equation}
Combining this with by Lemma \ref{estimates},
(\ref{positionvector1}), (\ref{monotone4}) and the fact that
$r>R_1$ in $\hbk\setminus D$,  we have
\begin{equation}\label{monotone7}
   \frac{d}{d\rho} \int_{\Sigma_\rho}(H_0-H)(W-W_0)\,d
\Sigma_\rho\ge 0.
\end{equation}
By Theorem \ref{positivemass}  and Corollary \ref{positivemass1},
we conclude that (\ref{positivity2}) is true.

Suppose equality holds in (\ref{positivity2}) for some future
directed null vector $\zeta$, then using Corollary
\ref{positivemass1}, we   have
\begin{equation}
\lim_{\rho\to\infty}\int_{\Sigma_\rho}(H_0-H)\mathbf{X}\cdot\zeta\,d\Sigma_\rho=0.
\end{equation}
By Proposition \ref{rigidity}, $(\Omega,g)$ is a domain in $\hbk$.
\end{proof}

\begin{rema}\label {monotone8} (\ref{monotone7}) means that for any
future directed null vector $\zeta$,
\begin{equation}\label{monotone9}
\frac{d}{d\rho}
\int_{\Sigma_\rho}(H_0-H)\mathbf{W}_{\Sigma_0}\cdot\zeta\,d
\Sigma_\rho\ge 0.
\end{equation}
From the proof of (\ref{monotone7}), it is easy to see that
(\ref{monotone9}) is still true if $\zeta$ is future directed
time-like.
\end{rema}
\begin{coro}\label{rigidity2} With the same assumptions and
notations as in Theorem \ref{positivity2}. Then for any
$\rho\ge0$,  the following vector is either zero or is future
directed non space-like:
 $$ \mathbf{m}(\rho)=\int_{\Sigma_\rho}(H_0-H)\mathbf{W}_{
\Sigma_0} \, d\Sigma_\rho.
$$
In particular, if $(\Omega,g)$ is not a domain in $\hbk$, then
$$
\mathbf{m}(\Omega)=\int_{\Sigma}(H_0-H)\mathbf{W}_{ \Sigma_0} \,
d\Sigma
$$
is future directed time-like.
\end{coro}
\begin{proof} By the proof of Theorem \ref{positivity1} and the
characterization of future directed non space-like vector
\cite{WY}, we conclude the first part of the corollary is true. If
$(\Omega,g)$ is not a domain in $\hbk$, by the rigidity part of
the theorem, this vector cannot be zero and cannot be null. Hence
it is future directed time-like.
\end{proof}

\begin{coro}\label{massmonotone} With the same assumptions and
notations as in Theorem \ref{positivity2}, let
$$\mathbf{m}(\rho)=\int_{\Sigma_\rho}(H_0-H)\mathbf{W}_{\Sigma_0}\,d\Sigma_\rho.
$$
then $\frac{d}{d\rho}\lf(|\mathbf{m}(\rho)|^2\ri)\geq 0$, where
$|\mathbf{m}(\rho)|$ is the Lorentz norm of $\mathbf{m}(\rho)$.
\end{coro}

\begin{proof} For any fixed $\rho_0$, let $\zeta$ be the vector
$\mathbf{m}(\rho_0)$. By Corollary \ref{rigidity2} as mentioned
above, $\zeta$ is a future directed non space-like, note that
\begin{equation}
\begin{split}
\frac{d}{d\rho} \lf(|\mathbf{m}(\rho)|^2\ri)\bigg|_{\rho=\rho_0}
=2\lf(\frac{d}{d\rho}\int_{\Sigma_\rho}(H_0 -H)\mathbf{W}_{\Sigma_0}
\, d\Sigma_\rho)\ri)\bigg|_{\rho=\rho_0}
\cdot\zeta\\
=2(\frac{d}{d\rho}\int_{\Sigma_t}(H_0
-H)\mathbf{W}_{\Sigma_0}\cdot\zeta d\Sigma_\rho)|_{\rho=\rho_0}
\end{split}\nonumber
\end{equation}
By   Remark \ref{monotone8}, we have
$$\frac{d}{d\rho}\lf(|\mathbf{m}(\rho)|^2\ri)\geq 0.
$$
\end{proof}

\begin{theo}\label{positivity3} Let $(\Omega, g)$ be as in Theorem
\ref{positivity1}  and let $\Sigma_0$ be the image of isometric
embedding of $\partial \Omega=\Sigma$ in
$\mathbf{H}^3_{-\kappa^2}$ which encloses $D$. Then for any $o\in
\Sigma_0$:

$$\int_{\Sigma_0}(H_0 -H)(y)\cosh\kappa r(o,y)\,d\Sigma_0(y) \ge 0
$$
Equality holds for some $o$ if and only if $(\Omega, g)$ is a
domain of $\mathbf{H}^3_{-\kappa^2}$. In particular, if $\Sigma$
is a standard sphere, then
$$
\int_{\Sigma_0}(H_0 -H)(y) \,d\Sigma_0(y) \ge 0
$$
and equality holds if and only if $(\Omega, g)$ is a domain of
$\mathbf{H}^3_{-\kappa^2}$.
\end{theo}
\begin{proof} We may embed $\hbk$ in $\Minkowski$ such that
$o=(0,0,0,1/\kappa)$. The result then follows from Theorem
\ref{positivity1}
\end{proof}
 Theorem \ref{positivity2} implies a previous result of the
 authors \cite{ST}
 \begin{theo}\label{Shi-Tam}
Let $(\Omega,g)$ be a compact three manifold with nonnegative
scalar curvature and with smooth boundary $\Sigma$. Suppose
$\Sigma$ has positive Gaussian curvature and positive mean
curvature $H$, then
$$
\int_\Sigma(H_0-H)\,d\Sigma\ge0
$$
where $H_0$ is the mean curvature of $\Sigma$ when it is
isometrically embedded in $\mathbb{R}^3$.
 \end{theo}

 This result follows from Theorem \ref{positivity3} and the
 following lemma.
 \begin{lemm} With the same assumptions and notations as in
 Theorem \ref{Shi-Tam}, suppose $H_{\kappa}$ is the mean curvature of $\Sigma$ when it
 is isometrically embedded in $\hbk$, $\kappa>0$. Then there exists $\kappa_i\to 0$ such that
 $\lim_{i\to\infty}H_{\kappa_i}=H_0.$
 \end{lemm}
\begin{proof} Let $\hbk$ be represented by the metric:
\begin{equation}\label{hbmetric}
\frac{4}{(1-\kappa^2|x|^2)^2}(dx_1^2+dx_2^2+dx_3^2)
\end{equation}
defined on $ x_1^2+ x_2^2+ x_3^2=|x|^2<\kappa^{-2}.$ We may assume
that $(\Sigma,g)$ is embedded to $\hbk$  with embedding
$\iota_\kappa$ such that $p$ is mapped to the origin, where $p$ is
some fixed point. Then it is easy to see that $\Sigma\subset
B_\kappa(0,2d)$ where $d$ the intrinsic diameter of $\Sigma$ and
$B_\kappa$ is the geodesic ball with respect to the metric
(\ref{hbmetric}).

Let $\iota_\kappa=(u_{1,\kappa},u_{2,\kappa},u_{3,\kappa})$ in
terms of the global coordinates $x_1,x_2,x_3$. Since the Gauss
curvature of $\Sigma$ is positive, by \cite[VI\S2,\S3]{P2}, we
conclude that for any $k\ge 0$, for $0<\kappa\le 1$ and $1\le j\le
3$, $|\nabla^ku_{j,\kappa}|$  are uniformly bounded. There exists
$\kappa_i\to0$ such that $\iota_{\kappa_i}$ together its
derivatives converge to an embedding of $\Sigma$ in the Euclidean
space. Using the fact that the embedding of $\Sigma$ in
$\mathbb{R}^3$ is unique, it is easy to see that the lemma is
true.
\end{proof}

Finally, we would like to give some examples to illustrate
that in certain situations, $\alpha$ in Theorem
\ref{positivity2} can be chosen as $1$. The proofs of these
examples are direct application of Theorem
\ref{positivity2},  Corollary \ref{positivity3} and the
representation of $\mathbf{X}$ given by
(\ref{positionvector1}).

\begin{exam}\label{example1} With the same assumptions and notations as in Theorem
\ref{positivity2}, if $\Sigma$ a standard sphere and $H$ is
constant, then   $m(\Omega, \zeta)\leq 0$ with
$\mathbf{W}_{\Sigma_0}=(x_1, x_2, x_3, t)$,   for all
future directed null vector $\zeta$. Here we have assumed
that $(0,0,0,1/\kappa)$ is inside $\Sigma_0$, which is the
image of $\Sigma$ under the embedding described in \S1.
\end{exam}

\begin{exam}\label{example2} With the same assumptions and notations as in Theorem
\ref{positivity2}, if $\Sigma$ a standard sphere   its mean
curvature $H$ is  orthogonal to the first eigenfunctions of
$\SS^2$, then $m(\Omega, \zeta)\leq 0$ with
$\mathbf{W}_{\Sigma_0}=(x_1, x_2, x_3, t)$, for all future
directed null vector $\zeta$. Here we assume that
$(0,0,0,1/\kappa)$ is the center of the geodesic ball in
$\hbk$ enclosed by $\Sigma_0$, which is the image of
$\Sigma$ under the embedding described in \S1.
\end{exam}

\begin{rema}It is easy to see that in Example \ref{example1} and
Example \ref{example2}, if $(\Omega,g)$ is not a domain of
$\mathbf{H}^3$, then on all of their small perturbations,
$m(\Omega, \zeta)\leq 0$ with $\mathbf{W}_{\Sigma_0}=(x_1,
x_2, x_3, t)$. For perturbations of  Example
\ref{example1}, we only assume that $(0,0,0,1/\kappa)$ is
inside $\Sigma_0$. For perturbations of Example
\ref{example2}, the inequality is true if
$(0,0,0,1/\kappa)$ is at some particular position inside
$\Sigma_0$.
\end{rema}

Let $f(p)=\int_{\Sigma_0}(H_0 -H)(y)\cosh\kappa
r(p,y)\,d\Sigma_0(y)$,  then it is a smooth function on
$D\subset\mathbf{H}^3_{-\kappa^2}$, and it would be interesting to
see some properties of this function. We first have:
\begin{prop}\label{criticalpoint}
Suppose $f$ has a critical point $o$ in the interior of $D$, and let
$o=(0,0,0,1/\kappa)$, then for $1\leq i \leq 3$,
$$\int_{\Sigma_0}(H_0 -H)\sinh\kappa r \cdot \phi_i =0.
$$
\end{prop}
\begin{proof}
It is easy to see that $\nabla_o r(o,y)$ is the unit tangential
vector of geodesic $\overrightarrow{yo}$, hence,
\begin{equation}
\begin{split}
\nabla_o r(o,y)&=(\cos\theta, \sin\theta \cos\psi, \sin\theta
\sin\psi, 0)\\
&=(\phi_1, \phi_2, \phi_3, 0)
\end{split}\nonumber
\end{equation}
Combine this fact and a direct computation, we see that the
conclusion is true.
\end{proof}
\begin{rema}Suppose $(\Omega, g)$ satisfying the assumptions in
Theorem \ref{positivity2}, and $o=(0,0,0,1/\kappa)$ be a critical
point of $f$, then Theorem \ref{positivity2} is true with
$\alpha=1$.
\end{rema}
Again by a direct computation we have
\begin{prop}Let $f$ be defined as above, then

$$\Delta f=3\kappa^2 f$$
\end{prop}
\begin{rema}
Suppose $(\Omega, g)$ satisfying the assumptions in Theorem
\ref{positivity2}, then by maximal principle, we know that $f$
cannot attain a local maximum inside of $D$; if there is $o\in D$
with $f=0$, then $f$ is identical to $0$ on the whole $D$ which
implies $(\Omega, g)$ is a domain of $\mathbf{H}^3_{-\kappa^2}$.
\end{rema}

\bibliographystyle{amsplain}

\end{document}